\documentstyle[12pt]{article}

\def\ad{a^{\dagger}}
\def\no{\hat{n}}
\title{Recursively minimally-deformed oscillators}
\author{
J. Katriel\thanks{
Permanent address: Department of Chemistry,
Technion, 32000 Haifa, Israel.
}
$\;$ and C. Quesne\thanks{
Directeur de recherches FNRS. email: cquesne@ulb.ac.be} \\
{\small \sl Physique Nucl\'eaire Th\'eorique et Physique
Math\'ematique} \\
{\small \sl Universit\'e Libre de Bruxelles} \\
{\small \sl Campus de la Plaine, CP229,}
{\small \sl B1050 Bruxelles, Belgium}}

\date{ }
\begin{document}
\setcounter{page}{1}

\maketitle

\pagestyle{plain}


















\begin{abstract}
A recursive deformation of the boson commutation relation is
introduced. Each step consists of a minimal deformation of
a commutator $[a,\ad]=f_k(\cdots;\no)$ into
$[a,\ad]_{q_{k+1}}=f_k(\cdots;\no)$, where $\cdots$ stands for the
set of deformation parameters that $f_k$ depends on, followed by
a transformation into the commutator
$[a,\ad]=f_{k+1}(\cdots,\, q_{k+1};\no)$ to which the deformed
commutator is equivalent within the Fock space. Starting from the
harmonic oscillator commutation relation $[a,\ad]=1$ we obtain the
Arik-Coon and the Macfarlane-Biedenharn oscillators at the
first and second steps, respectively, followed by a sequence
of multiparameter generalizations.
Several other types of deformed commutation relations related to
the treatment of integrable models and to parastatistics are also
obtained.  The ``generic'' form consists of a linear combination of
exponentials of the number operator, and the various recursive
families can be classified according to the number of free linear
parameters involved, that depends on the form of the initial
commutator.

\end{abstract}

\vspace{0.5cm}

\hspace*{0.3cm}
PACS: 0210, 0220, 0365

\newpage

\section{Introduction}

\hspace*{6 mm}
The study of deformed oscillators has already yielded
a plethora of formal results and applications, but the
attempts to introduce some order in the rich and
varied choice of deformed commutation (quommutation) relations
studied by different authors has so far achieved limited success.
Of particular interest in this respect are the treatments due to
Jannussis {\it et al.} \cite{Jan1,Jan2}, Daskaloyannis \cite{Das1,Das2},
McDermott and Solomon \cite{McD} and Meljanac {\it et al.}
\cite{Meljanac}.

The following is a partial list of deformations that have been
studied:
\begin{enumerate}
\item
The Arik-Coon oscillator \cite{Arik}
$$[a,\ad]_q\equiv a\ad -q\ad a=1$$
\item
The Macfarlane-Biedenharn oscillator \cite{Macfarlane,Biedenharn}
$$[a,\ad]_q=q^{-\no}$$
that has independently been proposed by Sun and Fu \cite{Sun}.
\item
The Chakrabarti-Jagannathan oscillator \cite{Chak}
$$[a,\ad]_p=q^{-\no}$$
\item
The Calogero-Vasiliev oscillator \cite{Vasiliev}
$$[a,\ad]=1+2\nu (-1)^{\no}\, ,$$
which for $2\nu=p-1$ is the Chaturvedi-Srinivasan parabose
oscillator of order $p$ \cite{Chaturvedi}.
\item
The Brzezi\'nski-Egusquiza-Macfarlane oscillator \cite{Brzezinski}
$$[a,\ad]=q^{-\no}(1+2\nu (-1)^{\no})$$
\item
Macfarlane's $q$-deformed Calogero-Vasiliev oscillator \cite{Mack}
$$a\ad-q^{\pm(1+2\nu K)}\ad a=[[1+2\nu K]] q^{\mp (\no+\nu-\nu K)}\, ,$$
where $K=(-1)^{\no}$ and $[[x]]=\frac{q^x-q^{-x}}{q-q^{-1}}$.
\end{enumerate}

In the present contribution we introduce a notion of minimal
deformation, that, along with the well known flexibility exhibited
by the presentation of the quommutation relations within the
Fock space, enables a
recursive deformation procedure to be formulated, generating
the various types of deformed oscillators listed above,
thus yielding a certain classification principle.
Moreover, the procedure proposed yields a multiparameter
generalization of the quommutation relations and
suggests that the ``generic'' structure involves
sums of exponentials of the number operator.

\section{Equivalence of quommutators and commutators}

\hspace*{6 mm}
Let $a$ and $\ad$ be two mutually conjugate operators and let $\no$
satisfy the commutation relations $[a,\no]=a$ and $[\no,\ad]=\ad$.
It follows that $\no$ commutes with $\ad a$ and with $a\ad$.
Furthermore, let
\begin{equation}
\label{comm}
a\alpha(\no)\ad -\ad \beta(\no) a=\gamma(\no)
\end{equation}
where $\alpha(\ell)$, $\beta(\ell)$ and $\gamma(\ell)$ are given
functions such that $\alpha(\ell)$ does not vanish for integral
and non-negative $\ell$.
This quommutation relation contains the form studied by McDermott and
Solomon~\cite{McD}, in which $\alpha(\hat n) = \gamma(\hat n) = 1$. It is a
symmetrized version of that studied by Meljanac {\it et al.}~\cite{Meljanac}
(corresponding to
$\alpha(\hat n) = 1$), to which it is easily shown to be equivalent. The
transformations introduced below take place within a Fock space representation
that is assumed to exist, possessing a non-degenerate ground state that
satisfies
$a|0>=0$. At least within this representation it is rather likely that the form
introduced by McDermott and Solomon \cite{McD} is sufficiently general. The
non-degeneracy requirement of the ground state has recently been relaxed by
several authors [16-18]
who introduced a doubly-degenerate ground state that
was found useful in the context of discussing intermediate
statistics. We shall not pursue this extension.
{}From the assumptions specified above it follows that
\begin{equation}
\ad |k>=\sqrt{F(k+1)}|k+1>
\end{equation}
and
\begin{equation}
a|k+1>=\sqrt{F(k+1)}|k>
\end{equation}
where
\begin{equation}
\label{FF}
F(k)=\sum_{i=0}^{k-1}
\frac{\gamma(i)\beta(i)\beta(i+1)\cdots\beta(k-1)}
{\alpha(i+1)\alpha(i+2)\cdots\alpha(k)\beta(k-1)} \; .
\end{equation}
Hence, $a\ad=F(\no +1)$, $\ad a=F(\no)$ and the
quommutator $[a,\ad]_Q\equiv a\ad-Q\ad a$ is
\begin{equation}
\label{five}
[a,\ad ]_Q= \frac{\gamma(\no)}{\alpha(\no +1)}+
\sum_{i=0}^{\no -1}
\frac{\gamma(i)\beta(i)\beta(i+1)\cdots\beta(\no -1)}
{\alpha(i+1)\alpha(i+2)\cdots\alpha(\no)}
\left(\frac{1}{\alpha(\no +1)}-\frac{Q}{\beta(\no -1)}\right)
\end{equation}
where $Q$ is arbitrary, but will usually be chosen to be equal
to unity, and where the appearance of the number operator within the
upper summation limit (as well as within the summand) has a well defined
meaning when applied to any Fock state.

Consider the following example:\hfill\break
Let $$a\ad-\ad q^{\no +1} a=1 \, ,$$
{\it i.e.,} $\alpha(\no)=1$, $\beta(\no)=q^{\no+1}$, $\gamma(\no)=1$.
This is equivalent to
$$[a,\ad ]=
1+(q^{\no}-1)\sum_{j=0}^{\no -1} q^{\frac{j(2\no -j-1)}{2}} \; .$$
The quommutation relation obtained for $q=-1$, {\it i.e.,}
$a\ad-\ad (-1)^{\no +1} a=1$, can be transformed with the aid of
the identity $$\sum_{j=0}^{2k} (-1)^{\frac{(j-1)j}{2}}=1$$
into the equivalent form $[a,\ad ]=(-1)^{\no}$, whose significance
was discussed by Quesne and Vansteenkiste \cite{QV}.

As a further example we consider the $q$-deformed Calogero-Vasiliev
oscillator, proposed by Macfarlane \cite{Mack}. This oscillator
can be transformed into
\begin{eqnarray}
\label{qCV}
[a,\ad]_Q &=& \frac{1}{2(q-q^{-1})}
\left\{ q^{\no}(q-Q)(q^{2\nu}+1)
+q^{-\no}(Q-q^{-1})(q^{-2\nu}+1) \right. \nonumber \\
 &+&\left. (-q)^{\no}(Q+q)(q^{2\nu}-1)
+(-q)^{-\no}(Q+q^{-1})(1-q^{-2\nu})\right\} \; .
\end{eqnarray}
This can be done either by starting from the quommutator
quoted in the introduction and applying the procedure illustrated
above, or, more simply, using the expressions for $a\ad$ and
for $\ad a$ presented by Macfarlane \cite{Mack}. In either case,
the expression obtained is written separately for $\no$ even and for
$\no$ odd, and the two expressions are combined with coefficients
of the form $\frac{1}{2}\Big(1+ (-1)^{\no}\Big)$
and $\frac{1}{2}\Big(1- (-1)^{\no}\Big)$, respectively.
Some further minor rearrangement yields eq. \ref{qCV},
that consists of a linear combination of four exponentials
in $\no$ (three, if $Q$ is chosen to be equal to $q$, $q^{-1}$,
$-q$ or $-q^{-1}$).

Another ``exotic'' quommutator is \cite{Solo1,Solo2}
$$a\ad-\frac{q^{\no+2}+1}{q(q^{\no}+1)}\ad a=1$$
{\it i.e.,} $\alpha(\no)=1$,
$\beta(\no)=\frac{q^{\no+3}+1}{q(q^{\no+1}+1)}$,
$\gamma(\no)=1$.
In this case
$$[a,\ad]_Q=1+\frac{q^{n+1}(q-Q)+1-qQ}{q^2-1}(1-q^{-\no})\, ,$$
which, for $q=Q$ reduces to the Macfarlane-Biedenharn oscillator
$[a,\ad]_q=q^{-\no}$.

\section{Recursive deformation of the harmonic oscillator}

\hspace*{6 mm}
In the following we will be interested in what appears to be a
somewhat more restricted framework.
Starting from $[a,\ad]=f_0(\no)$ let us assume that at the $k$th
step of a recursive procedure to be fully explicated below we
have obtained the commutation relation
$$[a,\ad ]=f_k(\no) \, .$$
We define the next minimal deformation to be
$$[a,\ad ]_{q_{k+1}}=f_k(\no) \, .$$
This minimally-deformed relation implies that, in the Fock-space
representation,
$$\ad |\ell>=\sqrt{F_{k+1}(\ell+1)}|\ell+1>$$
and $$a|\ell+1>=\sqrt{F_{k+1}(\ell+1)}|\ell>$$
where
\begin{equation}
\label{Fk}
F_{k+1}(\ell)=\sum_{i=0}^{\ell-1} q_{k+1}^{i}f_k(\ell-1-i) \, .
\end{equation}
It follows that
$$[a,\ad ]= f_{k+1}(\no)$$
where
$$f_{k+1}(\no) \equiv F_{k+1}(\no+1)-F_{k+1}(\no) \, .$$
This recurrence relation can also be written in the form
$$f_{k+1}(\no)=
\sum_{i=0}^{\no} q_{k+1}^{\no-i}\bigg(f_k(i)-f_k(i-1)\bigg)$$
provided that we define $f_k(-1) \equiv 0$. From the
recurrence relation it follows that if
$\lim_{q_1\rightarrow 1,q_2\rightarrow 1,\cdots,q_k\rightarrow 1}
f_k(\ell)=1$ for $\ell=0,\, 1,\, \cdots$, then
$\lim_{q_1\rightarrow 1,q_2\rightarrow 1,\cdots,q_{k+1}\rightarrow 1}
f_{k+1}(\ell) =1$. In other words, for all $k$, if $f_k(\no)$ is a
deformation of unity, so is $f_{k+1}(\no )$.

It will be convenient to define
$$\Phi_k(\ell)=
\sum_{{0\leq i_1,i_2,\cdots,i_k}\atop{(i_1+i_2+\cdots+i_k=\ell-k+1)}}
^{\phantom{{i_1}\atop{i_1}}}
q_1^{i_1}q_2^{i_2}\cdots q_k^{i_k} \, ,$$
which is easily shown to satisfy the limiting property
$$\lim_{q_1\rightarrow 1,q_2\rightarrow 1,\cdots,q_k\rightarrow 1}
\Phi_k(\ell)={\ell \choose {k-1}} \, .$$
Note that
$$\Phi_1(\ell)=q_1^{\ell}$$
$$\Phi_2(\ell)=\frac{q_1^{\ell}-q_2^{\ell}}{q_1-q_2}=
\frac{q_1^{\ell}}{q_1-q_2}+\frac{q_2^{\ell}}{q_2-q_1}$$
$$\Phi_3(\ell)= \frac{q_1^{\ell}}{(q_1-q_2)(q_1-q_3)}
+\frac{q_2^{\ell}}{(q_2-q_1)(q_2-q_3)}
+\frac{q_3^{\ell}}{(q_3-q_1)(q_3-q_2)}$$
or, in general,
$$\Phi_k(\ell)=
\sum_{i=1}^k \frac{q_i^{\ell}}{{\prod_{m=1}^k}^{\prime}(q_i-q_m)}$$
the prime indicating that $m\neq i$.

Let us now take
$$f_0(\ell)=\left\{\begin{array}{ccc}
1 & {\mbox{for}} & \ell\geq 0 \\
0 & {\mbox{for}} & \ell<0
\end{array}\right.$$
{\it i.e.,} start from the conventional harmonic oscillator
commutation relation, $[a,\ad]=1$.
With this initial value it can be shown that for $k\geq 1$
\begin{equation}
\label{fk}
f_k(\no)=\sum_{j=0}^{k-1}(-1)^{k-1-j}{{k-1}\choose j}
\Phi_k(\no+j)=
\sum_{i=1}^k \omega_{k,i} q_i^{\no}
\end{equation}
where
$\omega_{k,i}={\prod_{m=1}^k}^{\prime}
\left(\frac{q_i-1}{q_i-q_m}\right) \, . $
Applying the residue theorem to the function
$$f(z)=\frac{(z-1)^{\ell}}{\prod_{m=1}^k (z-q_m)}$$
we obtain
\begin{equation}
\label{iden}
\sum_{i=1}^k \frac{(q_i-1)^{\ell}}{{\prod_{m=1}^k}^{\prime}
(q_i-q_m)}=
\left\{ \begin{array}{cl}
1 &\;\;  \ell=k-1 \\
0 &\;\;  0\leq\ell<k-1
\end{array}\right. \; .
\end{equation}
The case $\ell=k-1$ yields
$$\sum_{i=1}^k \omega_{k,i}=1 \, ,$$
which clarifies the significance of eq. \ref{fk}, suggesting
that the coefficients $\omega_{k,i},\;  i=1,\, 2,\, \cdots,\, k,$
are the weights in an appropriate average.
Substituting eq. \ref{fk} in eq. \ref{Fk} we obtain, for $k\geq 1$,
$$F_{k+1}(\ell)=\sum_{i=1}^k \frac{(q_i-1)^{k-1}}
{{\prod_{m=1}^{k+1}}^{\prime}(q_i-q_m)}(q_i^{\ell}-q_{k+1}^{\ell})
=\sum_{i=1}^{k+1} \frac{(q_i-1)^{k-1}}
{{\prod_{m=1}^{k+1}}^{\prime} (q_i-q_m)} q_i^{\ell} \, ,$$
where use was made of the identity
$\sum_{i=1}^{k+1}\frac{(q_i-1)^{k-1}}
{{\prod_{m=1}^{k+1}}^{\prime}(q_i-q_m)}=0$ that corresponds
to $\ell=k-2$ in eq. \ref{iden}.
Using the latter identity once more we obtain the equivalent form
\begin{equation}
\label{nine}
F_k(\ell)=\sum_{i=1}^k\omega_{k,i} [\ell]_{q_i}
\end{equation}
where $[\ell]_{q_i}=\frac{q_i^{\ell}-1}{q_i-1}$ is the Jackson
$q_i$-(basic) integer. $F_k(\ell)$ is the weighted average of the
Jackson $q$-deformations of the integer $\ell$, in the $k$
different bases $q_1,\, q_2,\, \cdots,\, q_k$.
Thus, $F_1(\ell)=\frac{q_1^{\ell}-1}{q_1-1}$,
$F_2(\ell)=\frac{q_1^{\ell}-q_2^{\ell}}{q_1-q_2}=
\frac{q_1^{\ell}-1}{q_1-q_2}+\frac{q_2^{\ell}-1}{q_2-q_1}$, etc.

Using the Jackson $q$-derivative
${}_q D_x g(x)\equiv \frac{g(qx)-g(x)}{x(q-1)}$
we introduce the multi-parameter $q$-derivative
$${}_{q_1q_2\cdots q_k}D_x \equiv
\sum_{i=1}^k \omega_{k,i} \, {}_{q_i}D_x \; ,$$
which is a weighted average over the corresponding
Jackson $q_i$-derivatives. In particular, the
Macfarlane-Biedenharn $q$-derivative is a weighted average
over Jackson $q$-derivatives with respect to $q$ and $q^{-1}$,
{\it i.e.,} $${}_q {\overline{D}}_x=\omega_q\, {}_qD_x
+ \omega_{q^{-1}}\, {}_{q^{-1}}D_x$$
where $\omega_q=\frac{q-1}{q-q^{-1}}=
\frac{q^{1/2}}{q^{1/2}+q^{-1/2}}$ and
$\omega_{q^{-1}}=\frac{q^{-1/2}}{q^{1/2}+q^{-1/2}}$.
The multi-parameter $q$-derivative satisfies
$${}_{q_1q_2\cdots q_k}D_x x^{\ell}= F_k(\ell) x^{\ell-1}\, , $$
that enables the introduction of a corresponding $q$-exponential.

Thus, the minimal deformation of the conventional harmonic oscillator
$$[a,\ad ]=1 \, ,$$ is the relation
$$[a,\ad ]_{q_{1}}=1\, ,$$
which is due to Arik and Coon \cite{Arik}.
It is easily found that
$F_1(\ell)=[\ell]_{q_{1}}\equiv \frac{q_1^{\ell}-1}{q_1-1}$
and $f_1(\no)=q_1^{\no}$, {\it i.e.,} the Arik-Coon oscillator
is equivalent with
\begin{equation}
\label{one}
[a,\ad ]=q_1^{\no}\, .
\end{equation}
This equivalence had been pointed out by Kumari
{\it et al.} \cite{Kumari}.
Eq. \ref{one} suggests that the Arik-Coon oscillator
gets more and more classical, with increasing $\no$, for $q_1<1$,
and more and more quantal for $q_1>1$. In other, more picturesque
words, we have an ``energy dependent Planck's constant''.
This feature was discussed in refs. \cite{In,In1}, where it was
referred to as the Tamm-Dancoff cut-off.

Continuing, we consider the minimal deformation of eq. \ref{one},
{\it i.e.,} $[a,\ad ]_{q_2}=q_1^{\no}$.
This is the Chakrabarti-Jagannathan \cite{Chak}
two parameter oscillator,
which for $q_1=q_2^{-1}$ reduces to the Macfarlane-Biedenharn
\cite{Macfarlane,Biedenharn} oscillator.
The equivalent commutation relation is
\begin{equation}
\label{two}
[a,\ad ]=f_2(\no)
\end{equation}
where
$$f_2(\no)=\Phi_2(\no+1)-\Phi_2(\no) \, .$$
When $q_1=q_2^{-1}$ this expression reduces to the commutator
$[a,\ad ]=\frac{q_1^{(\no+\frac{1}{2})}+q_1^{-(\no+\frac{1}{2})}}
{q_1^{\frac{1}{2}}+q_1^{-\frac{1}{2}}}$,
that (with a slight change of notation) was noted by
Floreanini and Vinet \cite{Vinet}.

The equivalence between $[a,\ad]_{q_2}=q_1^{\no}$ and
$[a,\ad]=f_2(\no)$, and the fact that $f_2(\no)$ is symmetric in
$q_1$ and $q_2$, implies the well known equivalence of
$[a,\ad]_{q_1}=q_2^{\no}$ and $[a,\ad]_{q_2}=q_1^{\no}$.
It is perhaps appropriate to emphasize that the latter equivalence,
like the former, is only valid within the Fock space.

The minimal deformation of eq. \ref{two} yields
$[a,\ad]_{q_3}=f_2(\no)$, which can be written in the equivalent
commutator form
\begin{equation}
\label{three}
[a,\ad]=f_3(\no)
\end{equation}
where
$$f_3(\no)=\Phi_3(\no+2)-2\Phi_3(\no+1)+\Phi_3(\no) \, .$$

Continuing the recursion we note that since $f_k(\no)$
is a symmetric polynomial in
$q_1,\, q_2,\, \cdots ,\, q_k$ ({\it cf.} eq. \ref{fk}),
it follows that the $k$ relations
$$[a,\ad]_{q_i}=
f_{k-1}(q_1,q_2,\cdots,q_{i-1},q_{i+1},\cdots,q_k;\no)
\;\;\;\;\;\;\;\; i=1,\,2,\, \cdots,\, k$$
are all satisfied simultaneously with
$[a,\ad]=f_k(q_1,q_2,\cdots,q_k;\no)$.
Here, the dependence on the parameters is shown explicitly.

The present multiparameter deformation refers to a single
coordinate, unlike the multiparameter quantum groups associated
with the $n$-dimensional quantum space quommutation
relations [26-28].

At the $k$th step of the recursion we obtain a
commutator which is equal to a linear combination of $k$
exponentials of the number operator, with coefficients that
are fixed by the construction formulated.
We shall now consider a more general starting point,
involving a commutator that is equal to some polynomial
in the number operator. It will be found that once the number
of recursions exceeds the degree of the polynomial, the resulting
commutator is again equal to a sum of exponentials, but now with a
greater flexibility in the choice of the coefficients.

Taking
$$f_0(\ell)=\left\{\begin{array}{ccc}
\alpha_{0,1}\ell+\alpha_{0,0} & {\mbox{for}} & \ell\geq 0\\
0 & {\mbox{for}} & \ell<0
\end{array}\right. \, , $$
{\it i.e.,} $[a,\ad]=\alpha_{0,1}\no+\alpha_{0,0}$, we obtain
$$ f_1(\no)=\alpha_{1,0}+\alpha_{1,1} q_1^{\no}$$
where $\alpha_{1,0}=\frac{\alpha_{0,1}}{1-q_1}$ and
$\alpha_{1,1}=\alpha_{0,0}-\frac{\alpha_{0,1}}{1-q_1}$.
Thus,
\begin{equation}
\label{Vas}
[a,\ad]=1+2\nu q_1^{\no}
\end{equation}
is the first recursion of the relation $[a,\ad]=f_0(\no)$ with
$$f_0(\ell)=\left\{\begin{array}{ccc}
(1-q_1)\ell+(1+2\nu) & {\mbox{for}} & \ell\geq 0\\
0 & {\mbox{for}} & \ell<0
\end{array}\right. \;\;\;\; .$$
In particular, the Calogero-Vasiliev oscillator corresponds to
eq. \ref{Vas} with $q_1=-1$.\hfill\break
The second recursion yields
\begin{equation}
\label{Bz}
[a,\ad]=f_2(\no)=\alpha_{2,1} q_1^{\no}+\alpha_{2,2} q_2^{\no}
\end{equation}
where $\alpha_{2,1}=\alpha_{1,1} \frac{q_1-1}{q_1-q_2}$
and $\alpha_{2,2}=\alpha_{1,0}+\alpha_{1,1}\frac{q_2-1}{q_2-q_1}$.
For $q_1=q^{-1}$, $q_2=-q^{-1}$, $q_3=q$,
$\alpha_{2,1}=1$, $\alpha_{2,2}=2\nu$, the
minimal deformation of eq. \ref{Bz} is the
Brzezi\'nski-Egusquiza-Macfarlane oscillator \cite{Brzezinski}.

The third recursion yields
$$f_3(\no)=\alpha_{3,1}q_1^{\no}+\alpha_{3,2}q_2^{\no}
+\alpha_{3,3}q_3^{\no}$$
where
\begin{eqnarray*}
& & \alpha_{3,1}=\frac{\alpha_{2,1}}{q_1-q_3}(q_1-1) \\
& & \alpha_{3,2}=\frac{\alpha_{2,2}}{q_2-q_3}(q_2-1) \\
& & \alpha_{3,3}=
\frac{(q_3-q_2)\alpha_{2,1}
+(q_3-q_1)\alpha_{2,2}}{(q_3-q_1)(q_3-q_2)}(q_3-1) \, ,
\end{eqnarray*}
etc.

Starting from a quadratic expression in the number operator
$$f_0(\ell)=\left\{\begin{array}{ccc}
\alpha_{0,2}\ell^2+\alpha_{0,1}\ell+\alpha_{0,0} &
{\mbox{for}} & \ell\geq 0\\ 0 & {\mbox{for}} & \ell<0
\end{array}\right. $$
we obtain
$$f_1(\no)=\alpha_{1,1}q_1^{\no}+\alpha_{1,2}\no+\alpha_{1,3}$$
where
\begin{eqnarray*}
& & \alpha_{1,1}=
\frac{q_1(\alpha_{0,2}+\alpha_{0,1})
+(\alpha_{0,2}-\alpha_{0,1})}{(q_1-1)^2}+\alpha_{0,0} \\
& & \alpha_{1,2}=\frac{2\alpha_{0,2}}{1-q_1} \\
& & \alpha_{1,3}=
\frac{q_1(\alpha_{0,2}+\alpha_{0,1})
+(\alpha_{0,2}-\alpha_{0,1})}{(q_1-1)^2}
\end{eqnarray*}
and
$$f_2(\no)=\alpha_{2,1}q_1^{\no}+\alpha_{2,2}q_2^{\no}+\alpha_{2,3}$$
with appropriately defined coefficients. The next recursion yields
$$f_3(\no)=\alpha_{3,1}q_1^{\no}+\alpha_{3,2}q_2^{\no}
+\alpha_{3,3}q_3^{\no}\, ,$$
where the coefficients $\alpha_{3,1}$, $\alpha_{3,2}$ and
$\alpha_{3,3}$, that can be expressed in terms of
$\alpha_{0,0}$, $\alpha_{0,1}$ and $\alpha_{0,2}$,
can be chosen to agree with the coefficients of the
$q$-deformed Calogero-Vasiliev oscillator, eq. \ref{qCV},
provided that $Q$ is chosen to have one of the four values
for which eq. \ref{qCV} reduces to a sum of three exponentials,
say $Q=q$, and $q_1$, $q_2$ and $q_3$ are chosen to be
$q^{-1}$, $-q$ and $-q^{-1}$, respectively.

Thus, starting with $f_0(\no)$ that is a polynomial of degree
$k$ in $\no$ we obtain, upon performing the recursive minimal
deformation procedure, polynomials of decreasing degrees in $\no$
combined with linear combinations of exponentials in $\no$.
After $k$ steps we obtain just a linear combination of exponentials,
but the original $k$th degree polynomial allows a corresponding
number of coefficients in the linear combination to be chosen at will.

\section{Normal ordering relations and multi-parameter
deformed Stirling numbers}

\hspace*{6 mm}
To derive a normal-ordering formula for the pair of operators
$a$ and $\ad$ satisfying $[a,\ad]_{q}=f(\no)$ we first
use the identity
$$[AB,C]_{q_1q_2}=A[B,C]_{q_2}+q_2[A,C]_{q_1}B$$
to derive the relation
\begin{equation}
\label{commu}
[a^{\ell},\ad]_{q^{\ell}}=\{\ell(\no)\}a^{\ell-1}
\end{equation}
where
$\{\ell(\no)\}\equiv \sum_{i=0}^{\ell-1} q^{\ell-1-i}f(\no+i) \, .$
For $f(\no)=1$ we obtain
$\{\ell(\no)\}_{1}= [\ell]_{q_1}=\frac{q_1^{\ell}-1}{q_1-1}$,
that for $q_1=1$ is equal to $\ell$.
For $f(\no)=q_1^{\no}$ we have
$\{\ell(\no)\}_2=q_1^{\no}[[\ell]]_{q_1,q_2}$
where $[[\ell]]_{q_1,q_2}=\frac{q_1^{\ell}-q_2^{\ell}}{q_1-q_2}$.

Now, taking $q=q_{k+1}$ and $f(\no)=f_k(\no)$ (eq. \ref{fk})
we obtain, for $k\geq 1$,
$$\{\ell(\no)\}_{k+1}=\sum_{i=1}^k (q_i)^{\no}
\frac{q_i^{\ell}-q_{k+1}^{\ell}}{q_i-q_{k+1}} \omega_{k,i} \, .$$
Thus, $\{\ell(\no)\}_{k+1}$ is a multiparameter, operator-valued
deformation of the integer $\ell$.

Using eq. \ref{commu} we obtain that the coefficients in the normal
ordering formula
\begin{equation}
\label{qS}
(\ad a)^m=\sum_{\ell=1}^m (\ad)^{\ell} C_{m,\ell}(\no) a^{\ell}
\end{equation}
satisfy the initial condition $C_{1,1}(\no)=1$ and the
recurrence relation
$$C_{m+1,\ell}(\no)=q^{\ell-1} C_{m,\ell-1}(\no+1)
+\{\ell(\no)\}_k C_{m,\ell}(\no) \, .$$
In the appropriate limits this relation reduces to the Stirling,
$q$-Stirling and operator-valued $q$-Stirling coefficients,
{\it cf.} ref. \cite{Kibler}.

The normally-ordered form of an expression of the type $(\ad a)^m$
is not invariant with respect to the different equivalent
commutation and quommutation relations that the corresponding
pair of operators $a$ and $\ad$ satisfies.
Starting from the commutation relation $[a,\ad]=f_k(\no)$
we obtain $$[a^{\ell},\ad]={\overline{\{\ell(\no)\}}_k}$$
where
${\overline{\{\ell(\no)\}}_k}=\sum_{i=0}^{\ell-1} f_k(\no+i)
=\sum_{i=1}^k \omega_{k,i}[\ell]_{q_i}$
and $[\ell]_{q_i}=\frac{q_i^{\ell}-1}{q_i-1}$.
Hence, a normal ordering expansion of the form of eq. \ref{qS}
is obtained, with the coefficient satisfying the
recurrence relation
$${\overline{C}}_{m+1,\ell}(\no)=
{\overline{C}}_{m,\ell-1}(\no+1)
+{\overline{\{\ell(\no)\}}_k}{\overline{C}}_{m,\ell}(\no) \, .$$

Thus, the Arik-Coon quommutation relation
gives rise to the $q$-Stirling numbers as the coefficients
in the normally-ordered expansion, but if the (equivalent)
commutation relation $[a,\ad]=q^{\no}$ is used to effect
the normal ordering, the coefficients are operator-valued.
As an illustration consider $(\ad a)^2$, which, in terms of
the Arik-Coon quommutator is given by
$$(\ad a)^2=q(\ad)^2 a^2+\ad a \, ,$$
whereas in terms of the equivalent commutator becomes
$$(\ad a)^2=(\ad)^2 a^2+\ad q^{\no} a \, .$$
These two normally-ordered expressions are related to one another via
the identity $q^{\no}=(q-1)\ad a +1$.

\section{The inverse problem}

\hspace*{6 mm}
The following inverse problem may sometimes be of interest:
Given a commutator of some form,  can it be transformed into
a quommutator that, in some sense, is of simpler form?
To motivate this problem we recall that the normal ordering
problem for the Arik-Coon oscillator $[a,\ad]_q=1$ yields the
$q$-Stirling numbers as coefficients, whereas the equivalent
commutator relation, $[a,\ad]=q^{\no}$, yields a normal
ordering expansion with a new type of $\no$ dependent
(``operator valued'') $q$-Stirling numbers.
Given the latter commutator, we may wish to obtain the
equivalent quommutator that, in this case, yields a simpler
normal-ordering formula.

Thus, given $[a,\ad]=\phi(\no)$, where $\phi(0)=1$,
it can be shown straightforwardly that
$a\ad-\ad \beta(\no) a=1$, where
$$\beta(\no)=\frac{\sum_{i=1}^{\no+1}\phi(i)}
{\sum_{i=0}^{\no} \phi(i)} \, .$$
As an example we take $\phi(\no)=q^{\no}$ that yields
$\beta(\no)=q$, thus reproducing the Arik-Coon quommutator.

A somewhat different inverse problem involves the transformation of
$$[a,\ad]=\phi(\no )$$
into the equivalent form
$$[a,\ad ]_Q=\Phi(\no ) \, ,$$
choosing $Q$ so as to make $\Phi(\no)$ as simple as possible,
for a given $\phi(\no)$.
Since in the above quommutation relation $\alpha(\no)=\beta(\no)=1$,
we obtain
$$\Phi(\no) =\phi(\no)+(1-Q)\sum_{i=0}^{\no-1} \phi(i) \, ,$$
({\it cf.} eq. \ref{five}).
Thus, $\phi(\no)=q^{\no}$ yields
$$\Phi(\no)=q^{\no}\left(\frac{q-Q}{q-1}\right)+\frac{Q-1}{q-1} \; .$$
The ``best choice'' is very clear in this case, {\it i.e.},
$Q=q$, yielding $\Phi(\no)=1$.

Taking $$\phi(\no)=\alpha q_1^{\no}+\beta q_2^{\no}$$
we obtain
$$\Phi(\no)=\alpha q_1^{\no}\left(\frac{q_1-Q}{q_1-1}\right)
+\beta q_2^{\no}\left(\frac{q_2-Q}{q_2-1}\right)
+\alpha\frac{Q-1}{q_1-1}+\beta\frac{Q-1}{q_2-1}$$
In this case we have two equally good choices of $Q$,
{\it i.e.}, $Q=q_1$ and $Q=q_2$. The former yields
$$[a,\ad]_{q_1}=\Phi(\no)=
\beta\left(\frac{q_2-q_1}{q_2-1}\right) q_2^{\no}
+\left(\alpha+\beta\frac{q_1-1}{q_2-1}\right) \, .$$

For the special case $\phi(\no)=1+2\nu p^{\no}$ we obtain
$\sum_{i=0}^{\ell-1}\phi(i)=\ell+2\nu \frac{p^{\ell}-1}{p-1}$,
so, setting $Q=p$ it follows that $\Phi(\no)=1+2\nu+(1-p)\no$.
Hence, $[a,\ad]=1+2\nu p^{\no}$ is equivalent with
$[a,\ad]_p=1+2\nu+(1-p)\no$. Taking $p=-1$ we obtain that
$[a,\ad]=1+2\nu (-1)^{\no}$ is equivalent with
$\{a,\ad\}=1+2\nu+2\no$, {\it cf.} ref. \cite{Goodison}.
The latter is related to the realization of $osp(n/2,R)$ in terms
of parabosons, presented by Palev \cite{Palev}.

\section{Non-Fock space representations of the deformed commutation
relations}

\hspace*{6 mm}
The equivalence between quommutators and corresponding commutators,
presented in section 2, is a central ingredient of the recursive
minimal deformation procedure introduced in section 3. It was
noted in section 2 that the transformation proposed is  being
carried out within the Fock space representation.
Deformed oscillator algebras are known to have  additional,
non-Fock space, representations \cite{Rideau,Chai} that are
characterized by the existence of a Casimir operator with
non-trivial eigenvalues \cite{Kulish,Oh}. We stress that
there is no reason to expect these non-Fock space representations
to be the same for different ways of writing the commutation
relation that are equivalent within the Fock space.
While we do not wish to delve in a detailed analysis of these
non-Fock space representations for the different algebras
discussed, the following general observations indicate some of the
features to be expected.

The algebra $[a,\ad]=f_k(\no)=F_k(\no+1)-F_k(\no)$ has a
Casimir operator
$$C_k=F_k(\no)-\ad a \; .$$
This can be shown by noting that
$[C_k,\ad]=\Big(F_k(\no)-F_k(\no-1)\Big)\ad-\ad[a,\ad]=0$.
In the Fock-space representation a state $|0>$ exists, for
which the relation $a|0>=0$ is satisfied. Since, furthermore,
$F_k(\no)|0>=F_k(0)|0>$ and $F_k(0)\equiv 0$, it follows that
within this representation $C_k$ has eigenvalue $0$.
The non-Fock representations are characterised by non vanishing
eigenvalues of the Casimir operator.

The minimal deformation of the algebra just discussed,
$[a,\ad]_{q_{k+1}}=f_k(\no)$, has a Casimir operator as well,
{\it i.e,} ${\tilde{C}}_k=\mu_k(\no)-\nu_k(\no)\ad a$,
where $\mu_k(\no)$ and $\nu_k(\no)$ should be determined so as
to satisfy the condition $[{\tilde{C}}_k,\ad]=0$.
By adding a suitable constant one can set $\mu_k(0)=0$,
so that the Casimir operator vanishes on the Fock space
representation.
To determine $\mu_k(\no)$ and $\nu_k(\no)$ we note that
$$[{\tilde{C}}_k,\ad]=
\Big(\mu_k(\no)-\mu_k(\no-1)-\nu_k(\no)f_k(\no-1)\Big)\ad
+\Big(\nu_k(\no-1)-q_{k+1}\nu_k(\no)\Big)(\ad)^2a \; .$$
A sufficient condition for the vanishing of $[{\tilde{C}}_k,\ad]$
is
\begin{eqnarray}
\label{al}
& & \mu_k(\no)-\mu_k(\no-1)=\nu_k(\no)f_k(\no-1)\\
\label{be}
& & \nu_k(\no)=q_{k+1}^{-1}\nu_k(\no-1)
\end{eqnarray}
 From eq. \ref{be} we obtain $\nu_k(\no)=q_{k+1}^{-\no}$,
where the normalization $\nu_k(0)=1$ (which is consistent
with the choice $\mu_k(0)=0$ made above) was chosen.
Consequently, eq. \ref{al} becomes a recurrence relation for
$\mu_k(\no)$, {\it i.e.,}
$\mu_k(\no)=\mu_k(\no-1)+q_{k+1}^{-\no}f_k(\no-1)$.
This recurrence relation, along with the initial condition
$\mu_k(0)=0$, is satisfied by
$$\mu_k(\no)=\sum_{i=1}^{\no} q_{k+1}^{-i}f_k(i-1)\; .$$

For $f_k(\no)=\sum_{i=1}^k\omega_{k,i}q_i^{\no}$,
{\it cf.} eq. \ref{fk}, we obtain
$$\mu_k(\no)=
q_{k+1}^{-\no}\sum_{j=1}^{k+1}\omega_{k+1,j}[\no]_{q_j}
=q_{k+1}^{-\no}F_{k+1}(\no)$$
where the last equality follows from eq. \ref{nine}.
It follows that
$${\tilde{C}}_k=q_{k+1}^{-\no}\Big(F_{k+1}(\no)-\ad a\Big)=
q_{k+1}^{-\no}C_{k+1}\, .$$
The Casimir operators introduced can be used to investigate
the non-Fock space representations of the various deformed
oscillators presented, along the lines of refs. [32-35].
\vfill\eject

\section{Conclusions}

\hspace*{6 mm}
A recursive minimal-deformation of a commutator into a
quommutator, followed by a transformation of the resulting
quommutator into a new commutator, to which it is equivalent
within the corresponding Fock space, has been introduced.
The familiar deformed oscillators have been obtained at
appropriate steps of this recursive construction, along
with multiparameter generalizations that would be difficult to
guess otherwise. This recursive scheme provides a classification
of the existing deformed-oscillators.
The multi-parameter generalizations may appeal to investigators
who would like to use the deformed oscillator framework
in order to fit molecular or nuclear vibrational
spectra, and in similar contexts in which further flexibility
would be useful. To what extent they offer hints of further
fundamental developments remains to be seen.

\vspace{1cm}

\noindent
{\bf Acknowledgement} This research was carried out in the framework
of a project supported by the Fonds National de la Recherche
Scientifique, Belgique. We thank Dr. T. Brzezi\'nski for critically
reading the manuscript.

\end{document}